\documentstyle[11pt]{article}

\title{$^*$FORCING}

\author{Garvin
Melles\thanks{Would like to thank Ehud Hrushovski
for supporting him with funds from NSF Grant DMS 8959511} \\Hebrew
University of Jerusalem}

\newtheorem{theorem}{Theorem}[section]
\newtheorem{defi}[theorem]{Definition}
\newtheorem{lemma}[theorem]{Lemma}
\newtheorem{coro}[theorem]{Corollary}

\newcommand{\al}{\alpha}
\newcommand{\be}{\beta}
\newcommand{\ga}{\gamma}
\newcommand{\vp}{\varphi}
\newcommand{\sub}{\subseteq}
\newcommand{\proof}{{\sc proof} \hspace{0.1in}}

\begin{document}
\mathsurround=.1cm
\maketitle

\begin{abstract}
Let $M$ be a transitive model of $ZFC$ and let ${\bf B}$ be a
$M\!$-complete Boolean
algebra in $M.$ (Possibly a proper class.) We define a generalized notion
of forcing with such Boolean algebras, $^*$forcing. We prove that 
\begin{enumerate}
\item If ${\bf G}$ is a $^*$forcing complete ultrafilter on ${\bf
B},$ then $M[{\bf G}]\models ZFC.$
\item Let $H\sub M.$ If there is
a least transitive model $N$ such that $H\in M,$ $Ord^M=Ord^N,$ and
$N\models ZFC,$ then we denote $N$ by $M[H].$ We show that all models
of $ZFC$ of the form $M[H]$ are $^*$forcing extensions of $M.$ As an
immediate corollary we get that $L[0^{\#}]$ is a $^*$forcing extension
of $L.$
\end{enumerate}
\end{abstract}

\section{Introduction}

In this paper we introduce a generalization of forcing, $^*$forcing.
If $M$ is a transitive model of $ZFC$ and ${\bf B}$ is a $M\!$-complete
Boolean algebra in $M,$ then a $^*$forcing extension of $M$ is a
transitive class of the form $M[{\bf G}]$ where ${\bf G}$ is an
$M\!$-complete ultrafilter on ${\bf B}$. 
If ${\bf B}$ is a set Boolean algebra, then $^*$forcing and forcing
coincide, so the interesting case is when ${\bf B}$ is a proper class.
One indication of the importance of $^*$forcing is given by the fact that
$L[0^{\#}]$ is a $^*$forcing extension of $L.$ (In fact we show much
more, as in the abstract.) In what follows we review Boolean valued
models, and show how Boolean valued models on names for sets lead to
pre-$\!^*$forcing, $^*$forcing, and forcing with Boolean
values for formulas. We prove that if $M[{\bf G}]$ is a $^*$forcing
extension of $M$ and ${\bf G}$ is $^*$forcing complete, then $M[{\bf
G}]\models ZFC.$ In the final
section we show that if $M[H]$ is the smallest transitive model of
$ZFC$ such that $Ord^M=Ord^{M[H]}$ and $H\in M,$ then there is a
$^*$forcing complete ultrafilter ${\bf H}$ on the  Boolean
algebra ${\bf B}_{form}$ (a proper class in $M$) such that $M[H]=M[{\bf H}].$

\section{Boolean Valued Models}

\begin{defi}
A Boolean valued model $M(B)$(of the language of set theory) is a triple
$(U,I,E)$ where $U$ is the universe and $I$ and $E$ are binary
functions with values in some Boolean algebra. Keeping with tradition
we write $I(x,y)$ as $||x=y||$ and $E(x,y)$ as $||x\in y||.$ $I$ and
$E$ must satisfy the following laws:
\begin{enumerate}
\item $||x=x||=1$
\item $||x=y||=||y=x||$
\item $||x=y||\bullet||y=z||\leq ||x=z||$
\item $||x\in y||\bullet||v=x||\bullet||w=y||\leq ||v\in w||$
\end{enumerate}
For each formula  and $a_1,\ldots,a_n$ from
$U$ we define the Boolean value of  as follows:
\begin{enumerate}
\item $||\neg\psi(a_1,\ldots,a_n)||=-||\psi(a_1,\ldots,a_n)||$
\item $||\varphi\wedge\psi(a_1,\ldots,a_n)||=||\varphi(a_1,\ldots,a_n)||\bullet||\psi(a_1,\ldots,a_n)||$
\item $||\exists x\varphi(x,a_1,\ldots,a_n)||=\sum\limits_{a\in
M(B)}||\varphi(a,a_1,\ldots,a_n)||$
\end{enumerate}
If $M(B)$ is a Boolean valued model (or even if $||\ \ ||$ is defined
only for the atomic formulas) we can with any ultrafilter $F$ on $B$
construct a model $M(B)/F,$ as follows:
Define an equivalence relation over $M(B)$ by letting
$$x\equiv y\ \ \hbox{ iff }\ \ ||x=y||\in F$$
Define the relation $E$ on the equivalence classes by
$$[x]\ E\ [y]\ \ \hbox{ iff }\ \ ||x\in y||\in F$$
Let $M(B)/F=(M(B)/\equiv\ ,\ E).$
\end{defi}

\begin{lemma}\label{BVM}
If $M(B)$ is a Boolean valued model and $F$ is an ultrafilter on $B$
such that for every formula $\varphi=\exists x\psi(x)$ and $a_1,\ldots,a_n$ from the
universe of $M(B),$ if 
$$||\exists x\psi(x,a_1,\ldots,a_n)||\in F$$
imples there is an $a$ in the universe of $M(B)$ such that 
$$||\psi(a,a_1,\ldots,a_n)||\in F$$
then for every formula $\theta,$ and $b_1,\ldots,b_m$ from the
universe of $M(B),$
$$M(B)/F\models \theta([b_1],\ldots,[b_m])\ \hbox{ iff }\
||\theta(b_1,\ldots,b_m)||\in F$$
\end{lemma}
\proof Like the proof of Lemma 18.1 in [Jech].

\section{$^*$Forcing and $^*$Forcing Complete Ultrafilters}

\begin{defi}
If $M$ is a transitive model of $ZFC$ and ${\bf B}$ is a class Boolean
algebra in $M$ then we define the Boolean valued names for $M,\ M^{\bf B}$ as
$$\bigcup\limits_{\alpha\in Ord^M}M^{B_{\alpha}}$$
where $B_{\al}=\big\{b\in {\bf B}\mid rank(b)\leq \al\big\}.$
\end{defi}

\begin{defi}
Let $M$ be a transitive model of $ZFC$ and let ${\bf N}$ be a class
in $M.$ We call ${\bf N}$ in the rest of this definition the class of
names. Let $P$ be a class of unary relation symbols.
We define the class of formulas $Form(M,{\bf N},P)$ as the smallest
class ${\bf C}$ such that 
\begin{enumerate}
\item Each symbol in $P$ is in ${\bf C},$ together with each symbol of
$P$ with a name substituted for the variable
\item $x=y$ and $x\in y$ are in ${\bf C},$ together with substitutions
of variables in $x=y$ and $x\in y$ by names
\item If $\varphi(x_0,x_1,\ldots,x_n,A)$ is in ${\bf C}$ then so are
$$\neg\varphi(x_0,x_1,\ldots,x_n,A)$$
$$\exists x_0\varphi(x_0,x_1,\ldots,x_n,A)$$
\item Let $I$ be a set. If for each $i\in
I,\ \varphi(x_0,x_1,\ldots,x_n,A_i)\in {\bf C},$  
 then
$$\bigwedge\limits_{i\in I}\varphi(x_0,x_1,\ldots,x_n,A_i)\ \ \ \in{\bf
C}$$
$$\bigvee\limits_{i\in I}\varphi(x_0,x_1,\ldots,x_n,A_i)\ \ \ \in{\bf C}$$ 
\end{enumerate}
Let $Sent(M,{\bf N},P)$ be the elements of $Form(M,{\bf N},P)$
without free variables. Let $Sent(M,{\bf N},P)(q.f.)$ be the
sentences in $Sent(M,M^{\bf N},P)$ without quantifiers and similarly for 
$Form(M,{\bf N},P)(q.f.).$ 
\end{defi}

\begin{defi}
A pre-$^*$forcing notion is a quadruple $(M,M^{\bf B},
||\ \ ||^*,{\bf G})$  such that 
\begin{enumerate}
\item $M$ is a transitive model of $ZFC$
\item $M^{\bf B}$ is the class of Boolean valued names for $M$
\item $||\ \ ||^*$ is a function with domain the elements of
$Sent(M,M^{\bf B})(q.f.)$ such that 
\begin{enumerate}
\item $||x=y||^*=0$ if $x\neq y$ and $||x=y||^*=1$ if $x=y$
\item $||x\in y||^*=b$ if $x\in dom\,y$ and $y(x)=b$ and $0$ otherwise
\item If $I$ is a set and for each $i\in I,$ $A_i$ is a set of names, then
\begin{enumerate}
\item $||\neg\psi(A_0)||^*=-||\psi(A_0)||^*$
\item $||\bigvee\limits_{i\in I}\varphi(A_i)||^*=\sum\limits_{i\in
I}||\varphi(A_i)||^*$
\item $||\bigwedge\limits_{i\in I}\varphi(A_i)||^*=\prod\limits_{i\in
I}||\varphi(A_i)||^*$
\end{enumerate}
\end{enumerate}
\item ${\bf G}$ is a $M\!$-complete ultrafilter on ${\bf B},$ i.e., it
has the property that if $X\sub {\bf B}$ and 
$$\sum X \in {\bf G}$$
then there is an $x\in X$ such that $x\in {\bf G}.$
\end{enumerate}
We denote $(M^{\bf B},||\ \ ||^*)$ by $M({\bf B}).$ 
\end{defi}

\begin{theorem}
If $(M,M^{\bf B},
||\ \ ||^*,{\bf G})$ is a pre-$^*$forcing notion, 
then for every $\varphi(A)\in Sent(M,M^{{\bf B}})(q.f.),$
$$M({\bf B})/{\bf G}\models \varphi(A)\ \ \hbox{ iff }\ \
||\varphi(A)||^*\in {\bf G}$$
\end{theorem}
\proof By induction on the complexity of $\varphi$ and 
the definition of pre- $^*$forcing notion.

\vspace{.1in}

\noindent So we have a notion of a condition forcing a formula, at
least for quantifier free elements
of $Sent(M,M^{{\bf B}}).$
The disadvantage with pre-$\!^*$forcing is that the model
$M({\bf B})/{\bf G}$
is not in general a model of extensionality. But the function 
$||\ \ ||^*$ induces a function $||\ \ ||$ so that if we denote the 
$^*$Boolean valued model $(M,M^{\bf B},||
\ \ ||)$ by $M^{\bf B},$ then extensionality holds in $M^{\bf B}/{\bf
G}$. ($(M^{\bf B},||\ \ ||)$ is just
the usual notion of Boolean valued model for forcing.)

\begin{defi}
If $(M,M^{\bf B},||\ \ ||^*,{\bf G})$ is a pre-$^*$forcing notion then
we define by induction on $\Gamma(\rho(x),\rho(y))$ where $\rho(x)$ is
the rank of $x$ and $\Gamma$ is the canonical well ordering of
$Ord\times Ord,$
$$||x\in y||=\sum\limits_{t\in dom\,y}||x=t||\bullet||t\in y||^*$$
$$||x=y||=\prod\limits_{t\in dom\,x}(-||t\in x||^*+||t\in y||)\bullet
\prod\limits_{t\in dom\,y}(-||t\in y||^*+||t\in x||)$$
We call $(M,M^{\bf B},||\ \ ||,{\bf G})$ a $^*$forcing notion. 
\end{defi}

\begin{theorem}
If $(M,M^{\bf B},||\ \ ||^*,{\bf G})$ is a pre-$^*$forcing notion then
$||\ \ ||$ and $M^{\bf B}/{\bf G}$ with respect to $||\ \ ||$ are
defined and 
$$M^{\bf B}/{\bf G}\cong M[{\bf G}]$$
by the isomorphism $h[x]=i_{\bf G}(x)$ where 
$$M[{\bf G}]=\big\{i_{\bf G}(x)\mid x\in M^{\bf B}\big\}$$
with $\in$ on $M[{\bf G}]$ interpreted as the real $\in$ relation and $i_{\bf G}$ is defined in the
usual manner.   
\end{theorem}
\proof See lemmata 18.2, 18.3, 18.9 and exercise 18.6 of [Jech].

\begin{theorem}
(The forcing theorem for $^*$forcing) If 
$\varphi(A)\in Sent(M,M^{\bf B})(q.f.),$ then    
$$M[{\bf G}]\models \varphi(A)\ \ \hbox{ iff }\ \
||\varphi(A)||\in {\bf G}$$
\end{theorem}
\proof For every sentence $\vp(A)\in Sent(M,M^{\bf B})$ there is a
corresponding statement $\vp'(A')$  such that
$||\vp(A)||=||\vp'(A')||^*$ and 
$$M^{\bf B}/{\bf G}\models \vp(A)\ \leftrightarrow\ M({\bf B})/{\bf
G}\models \vp'(A')$$
Furthermore, $\vp(A)\in Sent(M,M^{\bf B})(q.f.)\ \leftrightarrow\
\vp'(A')\in Sent(M,M^{\bf B})(q.f.)$   
If $(M,M^{\bf B},||\ \ ||^*,{\bf G})$ is a pre-$^*$forcing notion and
if $\vp'(A')\in Sent(M,M^{\bf B})(q.f.)$  we know
that 
$$M({\bf B})/{\bf G}\models \vp'(A')\ \leftrightarrow\
||\vp'(A')||^*\in {\bf G}$$
This lifts to
$$M^{\bf B}/{\bf G}\models \vp(A)\ \leftrightarrow\ ||\vp(A)||\in {\bf G}$$
Since $M[{\bf G}]$ and $M^{\bf B}/{\bf G}$ are isomorphic we get the
statement in the theorem.

\begin{defi}
Let $(M,M^{\bf B},||\ \ ||,{\bf G})$ be a $^*$forcing notion. 
Suppose that for every first order formula
$\vp(x_0,x_1,\ldots,x_n)$ and set of names
$\dot{a}_1,\ldots,\dot{a}_n$ there is a set
$B_{\vp(x_0,\dot{a}_1,\ldots,\dot{a}_n)}$ of names such that for
every $C\supset B_{\vp(x_0,\dot{a}_1,\ldots,\dot{a}_n)},$
$$\sum\limits_{\dot{b}\in
B_{\vp(x_0,\dot{a}_1,\ldots,\dot{a}_n)}}||\vp(\dot{b},\dot{a}_1,\ldots,\dot{a}_n)||=
\sum\limits_{\dot{c}\in C}||\vp(\dot{c},\dot{a}_1,\ldots,\dot{a}_n)||$$
where the definition of $||\ \ ||$
is extended to all first order fomulas by
induction on their complexity  as follows:
\begin{enumerate}
\item $||\neg\vp(\dot{a}_1,\ldots,\dot{a}_n)||=-||\vp(\dot{a}_1,\ldots,\dot{a}_n)||$
\item $||\vp(\dot{a}_1,\ldots,\dot{a}_n)\ \wedge\
\psi(\dot{a}_1,\ldots,\dot{a}_n)||=
||\vp(\dot{a}_1,\ldots,\dot{a}_n)||\bullet ||\psi(\dot{a}_1,\ldots,\dot{a}_n)||$ 
\item $$||\exists
x_0\vp(x_0,\dot{a}_1,\ldots,\dot{a}_n)||=\sum\limits_{\dot{b}\in
B_{\vp(x_0,\dot{a}_1,\ldots,\dot{a}_n)}}||\vp(\dot{b},\dot{a}_1,\ldots,\dot{a}_n)||$$
\end{enumerate}
Then we say $(M,M^{\bf B},||\ \ ||,{\bf G})$ is a forcing notion with
Boolean values for formulas.
\end{defi}

\begin{theorem}
If $(M,M^{\bf B},||\ \ ||,{\bf G})$ is a forcing notion with Boolean values
for formulas, then
$$M[{\bf G}]\models ZFC-P$$
\end{theorem}
\proof First we show that the forcing theorem holds for all first
order formulas by induction on the complexity of the formulas, and then
the usual proof that a forcing extension (using a set Boolean algebra)
satisfies $ZFC$ for all the axioms except power set goes through.

\begin{defi}
Let $S$ be a class. A formula $\varphi(v_1,\ldots,v_n)$ is said to be a
formula of set theory enriched by $S$ if it is a first order formula 
built from the usual
atomic formulas of set theory together with atomic formulas of the
form $s(x)$ for $s\in S$ where the interpretation of $s(x)$ is $x\in s.$ 
\end{defi}

\begin{defi}
If $\al<Ord^M,$ let $M_{\al}^{\bf B}$ be the names of rank $\leq\al,$
and let $M_{\al}^{B_{\al}}$ be the elements of $M^{B_{\al}}$ of rank
$\leq\al.$ 
\end{defi}

\begin{defi}
Let $\al<Ord^{M}$ and let $(M,M^{\bf B},||\ \ ||,{\bf G})$ be a $^*$forcing
notion. Let $\vp(x_1,\ldots,x_n)$  be a formula and let
$\dot{a}_1,\ldots,\dot{a}_n$ be names in $M_{\al}^{B_{\al}}.$ By
induction on the complexity of $\vp$ we define 
$$||\vp(\dot{a}_1,\ldots,\dot{a}_n)||_{\al}$$
as follows:
\begin{enumerate}
\item
$||\vp(\dot{a}_1,\ldots,\dot{a}_n)||_{\al}=||\vp(\dot{a}_1,\ldots,\dot{a}_n)||$
if $\vp$ is quantifier free.
\item $||\vp(\dot{a}_1,\ldots,\dot{a}_n)\ \wedge\
\psi(\dot{a}_1,\ldots,\dot{a}_n)||_{\al}=$
$||\vp(\dot{a}_1,\ldots,\dot{a}_n)||_{\al}
\bullet||\psi(\dot{a}_1,\ldots,\dot{a}_n)||_{\al}$  
\item
$||\neg\vp(\dot{a}_1,\ldots,\dot{a}_n)||_{\al}=-||\vp(\dot{a}_1,\ldots,\dot{a}_n)||_{\al}$ 
\item $||\exists
x\vp(\dot{a}_1,\ldots,\dot{a}_n)||_{\al}=\sum\limits_{b\in
M_{\al}^{B_{\al}}}||\psi(b,\dot{a}_1,\ldots,\dot{a}_n)||_{\al}$ 
\end{enumerate}
\end{defi}

\begin{defi}
Let $\al<Ord^{M}$ and let $(M,M^{\bf B},||\ \ ||,{\bf G})$ be a $^*$forcing
notion. Let $\vp(x_1,\ldots,x_n)$ be a formula. We define by induction
on the complexity of $\vp$ the meaning of $\vp(x_1,\ldots,x_n)$
reflects in $M_{\al}^{B_{\al}}.$
\begin{enumerate}
\item If $\vp(x_1,\ldots,x_n)$ is quantifier free then
$\vp(x_1,\ldots,x_n)$ reflects in every $M_{\al}^{B_{\al}}$
\item $\vp(x_1,\ldots,x_n)$ reflects in $M_{\al}^{B_{\al}}$ iff
$\neg\vp(x_1,\ldots,x_n)$ reflects in $M_{\al}^{B_{\al}}$
\item $\vp(x_1,\ldots,x_n)\ \wedge\
\psi(x_1,\ldots,x_n)$ reflects in $M_{\al}^{B_{\al}}$ iff
$\vp(x_1,\ldots,x_n)$ and $\psi(x_1,\ldots,x_n)$ reflect in
$M_{\al}^{B_{\al}}$
\item If $\vp(x_1,\ldots,x_n)=\exists x_0\psi(x_0,x_1,\ldots,x_n)$
then $\vp(x_1,\ldots,x_n)$ reflects in $M_{\al}^{B_{\al}}$ iff
\begin{enumerate}
\item $\psi(x_0,x_1,\ldots,x_n)$ reflects in
$M_{\al}^{B_{\al}}$
\item $\forall\be>\al\exists\gamma>\be$ such that $\psi(x_0,x_1,\ldots,x_n)$ reflects in
$M_{\gamma}^{B_{\gamma}}$ and if $b\in M_{\gamma}^{B_{\ga}}$ and
$||\psi(b,\dot{a}_1,\ldots,\dot{a}_n)||_{\ga}\in {\bf G},$ for some
$\dot{a}_1,\ldots,\dot{a}_n$ from $M_{\al}^{B_{\al}}$ then there is an
$a\in M_{\al}^{B_{\al}}$ such that
$$||\psi(a,\dot{a}_1,\ldots,\dot{a}_n)||_{\al}\in {\bf G}$$ 
\end{enumerate}
\end{enumerate}
\end{defi}

\begin{defi}
If $(M,M^{\bf B},||\ \ ||,{\bf G})$ is a $^*$forcing notion then we say
that ${\bf G}$ is $^*$forcing complete if
\begin{enumerate}
\item For every 
$\al<Ord^M,$ there is an ordinal $\be$ such that if $\dot{a}$ is a
name with domain $M_{\al}^{B_{\al}},$ then for some name $\dot{b}$ in
$M_{\be}^{B_{\be}}$ 
$$||\dot{a}=\dot{b}||\in {\bf G}$$
\item For every formula of set theory $\vp(x_1,\ldots,x_n)$ and for every
$\al<Ord^M,$ there is an ordinal $\be>\al$ such that $\vp(x_1,\ldots,x_n)$
reflects in $M_{\be}^{B_{\be}}$
\end{enumerate}
\end{defi}

\begin{theorem}
If $(M,M^{\bf B},||\ \ ||,{\bf G})$ is a $^*$forcing notion and ${\bf
G}$ is $^*$forcing complete, then
$$M[{\bf G}]\models ZFC$$
\end{theorem}
\proof The only nontrivial axioms to prove are choice, union, power
set, replacement and separation. Separation follows from replacement
and from separation it is not hard to prove union. 

\noindent \underline{Replacement} First we prove by induction on the
complexity of $\vp$ that if $\vp$ reflects in $M_{\al}^{B_{\al}}$ then
$\forall\ \dot{a}_1,\ldots,\dot{a}_n\in M_{\al}^{B_{\al}}$
$$M[{\bf G}]\models \vp(i_{\bf G}(\dot{a}_1),\ldots,i_{\bf
G}(\dot{a}_n))\ \leftrightarrow\ i_{\bf G}[M_{\al}^{B_{\al}}]\models
\vp(i_{\bf G}(\dot{a}_1),\ldots,i_{\bf G}(\dot{a}_n))$$
$$\leftrightarrow\ ||\vp(\dot{a}_1,\ldots,\dot{a}_n)||_{\al}\in {\bf G}$$
(The second $\leftrightarrow$ follows from the definition of $||\ \
||_{\al}.$) If $\vp$ is quantifier free then it follows immediately.
The only nontrivial step is if $\vp=\exists x\psi(x).$ Suppose $M[{\bf
G}]\models \vp(i_{\bf G}(\dot{a}_1),\ldots,i_{\bf G}(\dot{a}_n)).$
Then for some $\beta<Ord^M$ there is a $\dot{b}\in M_{\be}^{B_{\be}}$
such that
$$M[{\bf G}]\models \psi(i_{\bf G}(\dot{b}),i_{\bf G}(\dot{a}_1),\ldots,i_{\bf
G}(\dot{a}_n))$$
Since ${\bf G}$ is $^*$forcing complete there is a $\ga>\be$ such that
$\psi(x_0,x_1,\ldots,x_n)$ reflects in $M_{\ga}^{B_{\ga}}$ which
implies by the induction hypothesis that 
$$i_{\bf G}[M_{\ga}^{B_{\ga}}]\models \psi(i_{\bf G}(\dot{b}),i_{\bf G}(\dot{a}_1),\ldots,i_{\bf
G}(\dot{a}_n))$$
and $||\psi(\dot{b},\dot{a}_1,\ldots,\dot{a}_n)||_{\ga}\in {\bf G}.$ 
Since ${\bf G}$ is $^*$forcing complete and $\vp(x_1,\ldots,x_n)$
reflects in $M_{\al}^{B_{\al}},$ $\psi(x_0,x_1,\ldots,x_n)$ reflects
in $M_{\al}^{B_{\al}}$ and for some $\dot{a}\in M_{\al}^{B_{\al}},$
$$||\psi(\dot{a},\dot{a}_1,\ldots,\dot{a}_n)||_{\al}\in {\bf G}$$
and $i_{\bf G}[M_{\al}^{B_{\al}}]\models \psi(i_{\bf
G}(\dot{a}),i_{\bf G}(\dot{a}_1),\ldots,i_{\bf G}(\dot{a}_n))$ which
implies 
$$||\vp(\dot{a}_1,\ldots,\dot{a}_n)||_{\al}\in {\bf G}$$
and $i_{\bf G}[M_{\al}^{B_{\al}}]\models 
\vp(i_{\bf G}(\dot{a}_1),\ldots,i_{\bf G}(\dot{a}_n)).$
Now if $||\vp(\dot{a}_1,\ldots,\dot{a}_n)||_{\al}\in {\bf G}$  then
we know that since ${\bf G}$ is $M\!$-complete that for some $a\in M_{\al}^{B_{\al}},$
$$||\psi(\dot{a},\dot{a}_1,\ldots,\dot{a}_n)||_{\al}\in {\bf G}$$
Now by the induction hypothesis we get 
$$M[{\bf G}]\models \psi(i_{\bf G}(\dot{a}),i_{\bf G}(\dot{a}_1),\ldots,i_{\bf
G}(\dot{a}_n))$$
i.e., $M[{\bf G}]\models \vp(i_{\bf G}(\dot{a}_1),\ldots,i_{\bf
G}(\dot{a}_n)).$
Now let $F(x)=y$ be a function defined by the formula $\vp(x,y,i_{\bf
G}(\dot{a}_1),\ldots,i_{\bf G}(\dot{a}_n))$ with the $\dot{a}_i$ from
$M_{\al}^{B_{\al}}.$ Let $X$ be a set and let $i_{\bf G}(\dot{X})=X.$
Let $\be>\al$ such that $\vp(x,y,\dot{a}_1,\ldots,\dot{a}_n)$ reflects in
$M_{\be}^{B_{\be}}$ and $\dot{X}\in M_{\be}^{B_{\be}}.$ Let $\dot{Y}$
be the name with domain $M_{\be}^{B_{\be}}$ such that $\forall\
\dot{a}\in M_{\be}^{B_{\be}},$
$$\dot{Y}(\dot{a})=\sum\limits_{x\in
dom\,\dot{X}}||\varphi(x,\dot{a},\dot{a}_1,\ldots,\dot{a}_n)
||_{\be}\bullet ||x\in\dot{X}||$$
Then $F[X]=i_{\bf G}(\dot{Y}).$

\noindent \underline{Power Set} By point 1. of the definition of
$^*$forcing complete and by separation we know that for each
$\al<Ord^M$ that $P(i_{\bf G}(M_{\al}^{B_{\al}}))$ exists. Now
applying separation and the fact that if $X\in M[{\bf G}],$ then for
some $\al<Ord^M, X\sub i_{\bf G}[M_{\al}^{B_{\al}}]$ we get the power
set axiom holds in $M[{\bf G}].$

\noindent \underline{Choice} If $X$ is a name a set in $M[{\bf G}]$
and $f$ is the name with domain $=\big\{(\check{x},x)^{\bf B}\mid x\in
dom\,X\big\}$ such that $f((\check{x},x)^{\bf B})$ is uniformly
$1,$ then $i_{\bf G}(f)$ is a function from $dom\,X$ onto $i_{\bf
G}(X).$ Since $dom\,X$ is well ordered (as it is in $M$) $i_{\bf G}(f)$ 
induces a well order on $i_{\bf G}(X).$

\begin{coro}
Let $T$ be a theory extending $ZFC.$ If $(M,M^{\bf B},||\ \ ||,{\bf
G})$ is a $^*$forcing notion and 
${\bf G}$ is $^*$forcing complete, then
$$M[{\bf G}]\models T$$
iff $\forall\vp\in T,\ \exists\al<Ord^M$ such that $M_{\al}^{B_{\al}}$
reflects $\vp$ and $||\vp||_{\alpha}\in {\bf G}.$ 
\end{coro}

\noindent We would like necessary and sufficient conditions that
$M[{\bf G}]\models ZFC$ but the proof does not seem to go in the other
direction unless ${\bf G}$ is a class in $M[{\bf G}].$ But if ${\bf
G}$ is a class in $M[{\bf G}]$ then the $\big\{i_{\bf
G}(M_{\al}^{B_{\al}})\mid \alpha<Ord^M\big\}$ form a cumulative
hierarchy in $M[{\bf G}]$ and we can prove:

\begin{theorem}
If $M[{\bf G}]$ is a $^*$forcing extension of $M,$ ${\bf G}$ is a
class in $M[{\bf G}]$ and $M[{\bf
G}]\models ZFC$ then
${\bf G}$ is $^*$forcing complete.
\end{theorem}
\proof If ${\bf G}$ is a class in $M[{\bf G}]$ then the $\big\{i_{\bf
G}(M_{\al}^{B_{\al}})\mid \alpha<Ord^M\big\}$ form a cumulative
hierarchy in $M[{\bf G}]$ so that for every formula $\vp$ and for
every ordinal $\al<Ord^M,$ there is a $\be>\al$ such that $\vp$
reflects (in the usual sense) in $i_{\bf G}(M_{\be}^{B_{\be}}).$ Now
we can show by induction on the complexity of $\vp$ that $\vp$
reflects in $M_{\be}^{B_{\be}}.$ Since the power set axiom holds in
$M[{\bf G}]$ we know that condition $1.$ in the definition of
$^*$forcing complete must also hold, otherwise for some $\al$ the
collection of subsets of $i_{\bf G}(M_{\al}^{B_{\al}})$ form a class
in $M[{\bf G}].$

\vspace{.1in}
 
\noindent If we consider models of set theory with an addition predicate $G(x)$
meant as an interpretation for ${\bf G}$ we can get necessary and
sufficient conditions that $M[{\bf G}]\models ZFC_{G(x)}$ with the
interpretation  of $G(x)$ being ${\bf G}$ where we
define $ZFC_{G(x)}$ as $ZFC$ except that we allow in the separation
and  replacement schemata the formulas to be enriched by $G(x).$ 
So for the rest of this section if we write $M[{\bf G}]$ we mean the
model in the language $(=,\in,\ G(x))$ with the interpretation of
$G(x)$ being ${\bf G}.$ If $(M,M^{\bf B},||\ \ ||,{\bf G})$ is a
$^*$forcing notion then we can extend the definition of $||\ \ ||$ so
that $(M,M^{\bf B},||\ \ ||)$ is a $^*$Boolean valued model in the
language $(=,\in,\ G(x))$ so that if $\dot{a}$ is a name then
$$i_{\bf G}(\dot{a})\in {\bf G}\ \ \hbox{ iff}\ \ M[{\bf G}]\models
G(i_{\bf G}(\dot{a})\ \ \hbox{ iff }\
\ ||G(\dot{a})||\in {\bf G}$$
If $b\in {\bf B},$ let 
$$G(\check{b})=b$$
Otherwise, if $\dot{a}$ is a name in $M^{\bf B}_{\al},$ let
$$||G(\dot{a})||=\sum\limits_{b\in
B_{\al}}||\dot{a}=\check{b}||\bullet||G(\check{b})||$$
Now we can extend the definition of $||\ \ ||$ to sentences in
$Sent(M,M^{\bf B},\{G(x)\})(q.f.)$ in the natural way.

\begin{defi}
If $(M,M^{\bf B},||\ \ ||,{\bf G})$ is a $^*$forcing notion then we say
that ${\bf G}$ is $^*$forcing complete relative to formulas enriched
by $G(x)$ if
\begin{enumerate}
\item For every 
$\al<Ord^M,$ there is an ordinal $\be$ such that if $\dot{a}$ is a
name with domain $M_{\al}^{B_{\al}},$ then for some name $\dot{b}$ in
$M_{\be}^{B_{\be}}$ 
$$||\dot{a}=\dot{b}||\in {\bf G}$$
\item For every formula of set theory enriched by $G(x),$ $\vp(x_1,\ldots,x_n)$ and for every
$\al<Ord^M,$ there is an ordinal $\be>\al$ such that $\vp(x_1,\ldots,x_n)$
reflects in $M_{\be}^{B_{\be}}$
\end{enumerate}
\end{defi}

\begin{theorem}
If $(M,M^{\bf B},||\ \ ||,{\bf G})$ is a $^*$forcing notion then
$M[{\bf G}]\models ZFC_{G(x)}\ $ iff ${\bf G}$ is $^*$forcing
complete relative to formulas of set theory enriched by $G(x).$
\end{theorem}
\proof Similar.

\section{$^*$Forcing and models of $ZFC$ of the form $M[H]$}

\begin{defi}
If $U$ is a transitive set, then $Def_A(U)$ is the set of all subsets of $U$ 
definable by formulas of set theory enriched by $A(x)$ with
quantifiers restricted to $U$
and parameters from $U.$ 
\end{defi}

\begin{defi}
If $M$ is a transitive model of $ZFC$ and $H\subset M,$ then if there
is a least transitive model $N$ of $ZFC$ such that $H\in M$ and
$Ord^M=Ord^N,$ then we denote $N$ by $M[H].$
\end{defi}

\begin{theorem}
Let $M$ be a  transitive model of $ZFC.$ If $H\sub M$  and $M[H]$ exists, then
$$M[H]=\bigcup\limits_{\alpha\in Ord^M}L_{\alpha}(M\cup\{H\})$$
where 
\begin{enumerate}
\item $L_0(M\cup \{H\})=\emptyset$
\item $L_{\alpha}(M\cup\{H\})=\bigcup\limits_{\beta<\alpha}L_{\beta}(M\cup\{H\})$ if $\alpha$ 
is a limit
\item $L_{\alpha+1}(M\cup\{H\})=\bigcup\limits_{A\in M\cup\{H\}}Def_A L_{\alpha}(M\cup\{H\})$
\end{enumerate} 

\end{theorem}
\proof See [Jech] exercise 15.13.

\vspace{.1in}

\noindent The above theorem induces a natural notion of ${\bf names}$
and of interpretation function $i_H.$

\begin{defi}
The class ${\bf names}$ is $\bigcup\limits_{\alpha\in
Ord^M}names_{\alpha}$ where we define 
$name_{\alpha}$ by induction on the ordinals in $M$ as 
follows: 
\begin{enumerate}
\item $names_0=\emptyset$
\item $names_{\alpha}=\bigcup\limits_{\beta<\alpha}names_{\beta}$ if
$\alpha$ is a limit 
\item $\tilde{y}\in names_{\alpha+1}$ iff $\tilde{y}$ is a function
with domain  $names_{\alpha}$ and
$\tilde{y}(\tilde{z})=\varphi^{\alpha}(\tilde{z},\tilde{z}_1,\ldots,\tilde{z}_n)$
where $\varphi(v_0,v_1,\ldots,v_n)$ is a formula of set theory
enriched by $M\cup\{H\}$ and
the $\tilde{z}_i$ are in $names_{\alpha}$  
\end{enumerate}
\end{defi}

\begin{defi}
If $M[H]$ exists we define by induction on $\alpha<Ord^M,$ $i_H$ as follows:
\begin{enumerate}
\item $i_H(\emptyset)=\emptyset$
\item If $\tilde{y}$ is a name in $names_{\alpha+1}$ with the form as
in the definition of ${\bf names}$ above, then 
$$i_H(\tilde{y})=\big\{i_H(\tilde{z})\mid i_H[names_{\alpha}]\models
\varphi(i_H(\tilde{z}),i_H(\tilde{z}_1),\ldots,i_H(\tilde{z}_n))\big\}$$ 
where $i_H[names_{\alpha}]=\big\{i_H(\tilde{z})\mid \tilde{z}\in
names_{\alpha}\big\}=L_{\alpha}(M\cup \{H\})$ 
\end{enumerate}
\end{defi}

\begin{defi}
${\bf B}_{form}$ is the class $Sent(M,{\bf names},M\cup\{H\})(q.f.)$
made into a $M\!$-complete Boolean algebra by equating sentences logically equivalent
in propositional logic. We have 
\begin{enumerate}
\item $0=[\emptyset=\emptyset\ \wedge\ \emptyset\neq\emptyset]$
\item $1=[\emptyset=\emptyset\ \vee\ \emptyset\neq\emptyset]$
\item $-[\varphi]=[\neg\varphi]$
\item If $\big\{A_i\mid i\in I\big\}$ is a set of sets of names, then 
\begin{enumerate}
\item $\prod\big\{[\varphi(A_i)]\mid i\in
I\big\}=[\bigwedge\limits_{i\in I}\varphi_i(A_i)]$
\item $\sum\big\{[\varphi(A_i)]\mid i\in
I\big\}=[\bigvee\limits_{i\in
I}\varphi(A_i)]$ 
\end{enumerate}
\end{enumerate}
\end{defi}

\begin{defi}
We define for every $\al<Ord^M$ a function  $Asn_{\al}$ and a class
function ${\bf I}$ such that
\begin{enumerate}
\item $dom\,Asn_{\al}=$ symbols of the form
$\varphi^{\alpha}(x_1,\ldots,x_m,\tilde{a}_1,\ldots,\tilde{a}_n)$ where
$\tilde{a}_1,\ldots,\tilde{a}_n$ are in $names_{\alpha}$ 
\item $range\,Asn_{\al}\sub Form(M,{\bf names},M\cup\{H\})(q.f.)$
\item $dom\,{\bf I}={\bf names}$
\item $range\,{\bf I}\sub M^{{\bf B}_{form}}$
\end{enumerate}
We define $Asn_{\al}$ by induction on the complexity of $\vp.$ 
\begin{enumerate}
\item If $\vp$ is quantifier free, then
$$Asn_{\al}(\vp^{\al})=[\vp]$$
\item $Asn_{\al}(\neg\vp^{\al})=-[Asn_{\al}(\vp^{\al})]$
\item $Asn_{\al}(\vp^{\al}\wedge\psi^{\al})=[Asn_{\al}(\vp^{\al})]\bullet [Asn_{\al}(\psi^{\al})]$
\item $Asn_{\al}(\exists x\psi(x)^{\al})=[\bigvee\limits_{\tilde{a}\in
names_{\al}}Asn_{\al}(\psi^{\al}(\tilde{a}))]$
\end{enumerate}
We let ${\bf
I}=\bigcup\limits_{\al\in Ord^M}I_{\al}$ where we define 
$I_{\al}$ by induction on $\al$ as follows: Let
$I_0(\emptyset)=\emptyset.$ Suppose $I_{\alpha}$ has been defined. 
Now if $\tilde{y}\in name_{\al+1}$ with elements in the range of
$\tilde{y}$ of the form
$\vp^{\alpha}(\tilde{a},\tilde{a}_1,\ldots,\tilde{a}_n)$ then
$I_{\al+1}(\tilde{y})$ has
domain
$$\big\{I_{\alpha}(\tilde{a})\mid \tilde{a}\in names_{\al}\big\}$$
and
$$I_{\al+1}(\tilde{y})(I_{\al}(\tilde{a}))=
Asn_{\al}(\vp(\tilde{a},\tilde{a}_1,\ldots,\tilde{a}_n))$$    
\end{defi}

\begin{theorem}
Every model of the form $M[H]$ is a $^*$forcing extension of $M$ via
${\bf B}_{form}$ by a $^*$forcing complete ultrafilter. 
\end{theorem}
\proof Let ${\bf H}=\big\{[\vp(A)]\mid \vp(A)\in Sent(M,{\bf
names},M\cup\{H\})(q.f.)\ \wedge\ M[H]\models \vp(i_{H}[A])\big\}.$
We leave to the reader the proof that ${\bf H}$ is an ultrafilter and 
${\bf H}$ is $M\!$-complete. 
Since ${\bf H}$ is definable
inside $M[H]$ (If $\vp(A)\in Sent(M,{\bf
names},M\cup\{H\})(q.f.)$ then  $M[H]\models \vp(i_{H}[A])\big\}$ iff
$\exists\al(A\sub names_{\al}\ \wedge\ i_H[names_{\al}]\models \vp(i_{H}[A])$)
we know that $M[{\bf H}]=i_{\bf H}[M^{{\bf
B}_{form}}]\sub M[H].$  
To prove the other direction we first prove by induction on
the rank of $\tilde{y}$ that 
$$i_H(\tilde{y})=i_{\bf H}({\bf I}(\tilde{y}))$$
For $\tilde{y}=\emptyset$ it follows from the definitions. So let
$\al=\be+1$ and let $\tilde{y}\in names_{\al}-names_{\be}.$ Suppose
$\tilde{y}(\tilde{a})=\vp^{\al}(\tilde{a},\tilde{a}_1,\ldots,\tilde{a}_n).$
Then
$$L_{\al}(M\cup\{H\})\models\vp(i_H(\tilde{a}),i_H(\tilde{a}_1),\ldots,i_H(\tilde{a}_n))\
\hbox{ iff}$$
$$L_{\al}(M\cup\{H\})\models Asn_{\al}\vp(i_H(\tilde{a}),i_H(\tilde{a}_1),\ldots,i_H(\tilde{a}_n))\
\hbox{ iff}$$
$$M[H]\models Asn_{\al}\vp(i_H(\tilde{a}),i_H(\tilde{a}_1),\ldots,i_H(\tilde{a}_n))\
\hbox{ iff}$$
$$Asn_{\al}\vp(\tilde{a},\tilde{a}_1,\ldots,\tilde{a}_n)\in {\bf H}$$
As a result,
$$i_H(\tilde{y})=\big\{i_H(\tilde{a})\mid
L_{\al}(M\cup\{H\})\models\vp(i_H(\tilde{a}),i_H(\tilde{a}_1),\ldots,i_H(\tilde{a}_n))\big\}$$
$$=\big\{i_H(\tilde{a})\mid
Asn_{\al}\vp(\tilde{a},\tilde{a}_1,\ldots,\tilde{a}_n)\in {\bf H}\big\}$$
$$=\big\{i_{\bf H}({\bf I}(\tilde{a}))\mid
Asn_{\al}\vp(\tilde{a},\tilde{a}_1,\ldots,\tilde{a}_n)\in {\bf H}\big\}$$
$$=i_{\bf H}({\bf I}(\tilde{y}))$$
So we have $M[H]=i_H[{\bf names}]\sub i_{\bf H}[M^{{\bf
B}_{form}}]=M[{\bf H}].$ (As $M[H]\models ZFC$ and ${\bf H}$ is a
class in $M[H]$ we know ${\bf H}$ must
be $^*$forcing complete.)

\vspace{.1in}

\noindent As a corollary, we get the following.

\begin{coro}
$L[0^{\#}]=L[{\bf G}]$ where ${\bf G}$ is a $^*$forcing complete
ultrafilter on ${\bf B}_{form}^L.$ 
\end{coro}

\vspace{.2in}

\begin{center}
REFERENCES
\end{center}

\noindent [Jech] T. Jech, {\em Set Theory}, Academic Press, 1978.

\end{document}